\documentclass{amsart}

 \usepackage{geometry}
\usepackage{amssymb}
\usepackage{amsthm,amscd,amsmath}
\usepackage{tikz-cd}
\usepackage{color}
\usetikzlibrary{decorations.pathmorphing}

\usepackage{enumerate}
\usepackage{cite}
\usepackage{mathrsfs}
\usepackage[pagebackref]{hyperref}

\newcommand{\bpm}{\begin{pmatrix}}
\newcommand{\epm}{\end{pmatrix}}

\newcommand{\D}{\mathscr{D}}
\newcommand{\cat}[1]{\mathbf{#1}}
\newcommand{\twocat}[1]{\underline{\mathbf{#1}}}

\newcommand{\scat}[1]{\mathbf{#1}_\bullet}
\newcommand{\mc}{\mathcal}

\newcommand{\del}{\partial}
\newcommand{\iso}{\cong}
\renewcommand{\equiv}{\simeq}
\renewcommand{\subset}{\subseteq}

\newcommand{\map}[1]{\stackrel{ {#1} }{ \longrightarrow }}

\newcommand{\Ho}{\mathrm{Ho}}
\newcommand{\HO}{\mathrm{HO}}
\newcommand{\ob}{\mathrm{ob}}

\newcommand{\id}{\mathrm{id}}
\newcommand{\op}{\mathrm{op}}

\newcommand{\dia}{\mathrm{dia}}
\newcommand{\colim}{\operatornamewithlimits{colim}}
\newcommand{\mono}{\rightarrowtail}

\newcommand{\ladj}{\dashv}
\newcommand{\diag}{\mathrm{diag}}
\newcommand{\eq}{\mathrm{eq}}
\newcommand{\str}{\mathrm{str}}

\newcommand{\ddown}[1]{\Delta\downarrow{#1}}

\renewcommand{\a}{\alpha}
\renewcommand{\b}{\beta}

\newcommand{\g}{\gamma}

\newcommand{\s}{\sigma}
\newcommand{\ph}{\varphi}

\newtheorem{thm}[equation]{Theorem}
\newtheorem{cor}[equation]{Corollary}
\newtheorem{lemma}[equation]{Lemma}
\newtheorem{prop}[equation]{Proposition}

\newtheorem{conj}[equation]{Conjecture}
\theoremstyle{definition}
\newtheorem{definition}[equation]{Definition}

\theoremstyle{remark}
\newtheorem{remark}[equation]{Remark}

\title{A higher Whitehead theorem and the embedding of quasicategories in prederivators}
\author{Kevin Arlin\textsuperscript{\scriptsize 1}}
\thanks{\textsuperscript{\scriptsize 1}n\'e Carlson}
\email{kdcarlson@ucla.edu}
\address{Mathematics Department,
         University of California, Los Angeles,
         Box 951555,
         Los Angeles, CA 90095,
         USA}
\begin{document}

\begin{abstract}
We prove a categorified Whitehead theorem showing that the 2-functor $\mathrm{HO}$
associating a prederivator to a quasicategory reflects equivalences.
The question of whether $\mathrm{HO}$ is bicategorically fully faithful (that is, whether morphisms
and 2-morphisms can be uniquely lifted from prederivators to quasicategories)
is more subtle. We can show that small quasicategories embed fully faithfully,
both bicategorically and with respect to a certain simplicial enrichment, into 
prederivators defined on arbitrary small categories. When the
quasicategories are not necessarily small, or when the prederivators are
defined only on homotopically finite categories, the 2-categorical argument breaks down,
although the simplicial version continues to go through. We give a conjectural 
counterexample to bicategorical full faithfulness in general. \end{abstract}


\maketitle

\section{Introduction}

Every notion of an abstract homotopy theory $\mathcal C$, whether
an $\infty$-category or a model category, admits the common 
underlying structure the homotopy category $\Ho(\mc C)$. Moreover, for every 
category~$J$ there exists a homotopy theory $\mc C^J$ of $J$-shaped diagrams 
in $\mc C$, which thus has its own homotopy category $\Ho(\mc C^J)$.
Indeed, each homotopy theory $\mc C$ gives rise to a 2-functor 
$\Ho(\mc C^{(-)})$ sending categories to categories. 
This is known as the ``prederivator" of $\mc C$. (Pre)derivators were, in fact,
axiomatized independently by Grothendieck \cite{grothendieck}, Heller 
\cite{heller} and Franke \cite{Franke}, before the modern development of
flexible models of $\infty$-categories.

Prederivators are thus often treated in the literature as a notion of
abstract homotopy theory, but this intuition has not always been
referred to mathematical fact. One might naturally ask for an embedding of 
quasicategories in prederivators, in any of various homotopical senses. 
It is clear that not every prederivator is levelwise equivalent to the prederivator
associated to a quasicategory. Fuentes-Keuthan, Kedziorek, and Rovelli  have recently described the image
up to isomorphism of prederivators \cite{rovelli}, but there is as yet no proposed description 
of the image up to equivalence.
We view the latter problem as the key remaining question in this area, 
which will be investigated in future work on 2-categorical Brown 
representability.

The question of embedding quasicategories into prederivators bifurcates
into two approaches, as quasicategories naturally form a simplicially 
enriched category $\scat{QCat}$, while prederivators naturally form only a 2-category 
$\twocat{PDer}$. We compare both $\scat{QCat}$ to the simplicially 
enriched category $\scat{PDer}$ of \cite{muro} and $\twocat{PDer}$ to the 2-category
$\twocat{QCat}$ of quasicategories first studied by Joyal. 

The reason for investigating the 2-categorical as well as the simplicial
comparison is that the 2-category $\twocat{PDer}$ is to understand to
what extent the theory of quasicategories can be reduced to ordinary
category theory: $\twocat{PDer}$ is an ordinary 2-category of 2-functors
into $\twocat{Cat}$.

\subsection*{Summary of Results}

We will denote henceforth by $\twocat{QCAT}$ the 2-category of quasicategories. A prederivator 
is a 2-functor $\D:\twocat{Dia}^\op\to\twocat{CAT}$, where the indexing category $\twocat{Dia}$
is generally either the 2-category $\twocat{Cat}$ of small categories or that 
$\twocat{HFin}$ of homotopy finite categories. Denoting the 2-category of prederivators by 
$\twocat{PDER}$, we can construct 2-functors $\HO:\twocat{QCAT}\to\twocat{PDER}$ for each 
$\twocat{Dia}$.

Our results are as follows.

\emph{First result:} In Theorem \ref{Maintheorem}, we show that
the 1-category $\cat{QCAT}$ of quasicategories embeds fully 
faithfully in any category $\cat{PDER}^\str$ of prederivators with strictly 2-natural morphisms. 
This extends to an embedding $\scat{QCAT}\to\scat{PDER}$ of simplicial categories,
where the domain has the usual simplicial enrichment. Thus, 
quasicategories and their mapping spaces
 can be recovered  \emph{up to isomorphism} from their prederivators and strict 
 maps. The on-the-nose quality of this statement 
 reflects the use of strict transformations of prederivators, as opposed to the pseudo-natural
 transformations that appear in the later results.

\emph{Second result:} In Theorem \ref{Mainthm2cat}, we show that 
$\HO:\twocat{QCAT}\to\twocat{PDER}$ is bicategorically fully faithful when
restricted to small quasicategories, as long as $\twocat{Dia}$ contains all of 
$\twocat{Cat}$.
The main tool is Joyal's delocalization theorem, as published by Stevenson \cite{stevenson}, 
which we recall as Theorem \ref{delocalization}. This allows us to write every quasicategory
as a localization of a 1-category. 

This result is similar to the theorem of Renaudin \cite{Ren} that a certain
2-category of combinatorial model categories embeds bicategorically fully faithfully in
$\twocat{PDER}_!$, the 2-category of cocomplete prederivators with cocontinuous
morphisms, insofar as combinatorial model categories model locally presentable
quasicategories, cocontinuous maps out of which are determined by restriction
to small dense subcategories. 

\emph{Third  result:} Our final result, applicable in more generality, is Theorem 
\ref{whiteheadforquasicats}, which shows that every version of 
$\HO$ is bicategorically \emph{conservative}, in that equivalences of 
quasicategories are reflected by $\HO$. 
In other words, the prederivator is enough to distinguish equivalence classes of 
abstract homotopy theories, no matter which size choices we make. The proof is 
unrelated to that of Theorem \ref{Mainthm2cat}, and relies on the author's 
Whitehead theorem for the 2-category of unpointed spaces
\cite{whiteheadforspaces}.

\emph{Conventions:} We will denote the category, the 2-category,
and the simplicial category of foos respectively by
\[\cat{foo},\twocat{foo},\scat{foo}\]
Furthermore, when applicable, the above will designate
the category of \emph{small} foos while
\[\cat{FOO},\twocat{FOO},\scat{FOO}\]
will refer to \emph{large} ones. We operationalize the term \emph{large}
to mean ``small with respect to the second-smallest Grothendieck
universe." 

For us 
2-categories are strict: they have strictly associative composition
and strict units preserved on the nose by 2-functors. We denote the horizontal composition of 2-
morphisms by $*$, so that if $\a:f\Rightarrow g:x\to y$ and
$\b:h\Rightarrow k:y\to z$, we have $\b*\a:h\circ f\Rightarrow k\circ g$.
If $\mc C$ is a category (or a 2-category, simplicially enriched category, etc)
with objects $c_1$ and $c_2,$ we denote the set (or category, simplicial set, etc)
of morphisms by~$\mc C(c_1,c_2)$.

\emph{Acknowledgements:} Thanks to an anonymous reviewer, whose detailed comments on an earlier
draft contributed greatly to whatever readability this manuscript may possess;
Paul Balmer, for his advice and support; Denis-Charles Cisinski, for suggesting the application 
of the delocalization theorem; James Richardson, for pointing out an error in the original proof 
of Theorem \ref{whiteheadforquasicats}, to an anonymous referee for exceptionally detailed comments 
which greatly improved the paper's readability, and to Martin Gallauer and Mike Shulman for helpful 
comments and conversations.

\section{Background on 2-categories and prederivators}

Below we recall the various 2-categorical definitions we will require.

\begin{definition}\label{2Cat}
Morphisms between 2-functors will be either 2-natural or pseudonatural transformations depending on context.
Let us recall that, if $\mc K,\mc L$ are 2-categories and $F,G:\mc K\to\mc L$ are 2-functors,
a pseudonatural transformation $\Lambda:F\Rightarrow G$ consists of 
\begin{itemize}
\item Morphisms $\Lambda_x:F(x)\to G(x)$ associated to every object $x\in\mc K$
\item 2-morphisms $\Lambda_f: \Lambda_y\circ F(f)\Rightarrow G(f)\circ 	\Lambda_x$ for every morphism
$f:x\to y$ in $\mc K$
\end{itemize}
satisfying the coherence conditions
\begin{itemize}
\item (Pseudonaturality) $\Lambda_f$ is an isomorphism, for every $f$.
\item (Coherence) $\Lambda$ is a functor from the underlying 1-category of $\mc K$ to the category
of pseudo-commutative squares in $\mc L$, that is, squares commuting up to a chosen isomorphism, where composition is by pasting. 
\item (Respect for 2-morphisms) For every 2-morphism $\a:f\Rightarrow g:x\to y$ in $\mc K$, we have the equality of 2-morphisms
\[\Lambda_g\circ(\Lambda_y*F(\alpha))=(G(\alpha)*\Lambda_x)\circ\Lambda_f:\Lambda_y\circ F(f)\Rightarrow G(g)\circ \Lambda_x.\]
\end{itemize}
In case all the $\Lambda_f$ are identities, we say that $\Lambda$ is strictly 2-natural, in
which case the axiom of coherence is redundant, and that of respect for 2-morphisms becomes simply
$\Lambda_y*F(\alpha)=G(\alpha)*\Lambda_x$.
\end{definition}

\begin{definition}\label{mod}
The morphisms between pseudonatural transformations, are called \emph{modifications.} 
A modification 
$\Xi:\Lambda\Rrightarrow \Gamma:F\Rightarrow G:\mc K\to \mc L$ 
consists of 2-morphisms $\Xi_x:\Lambda_x\to\Gamma_x$ for each object $x\in \mc K$, subject to 
the condition
\[(G(f)*\Xi_x)\circ\Lambda_f=\Gamma_f\circ(\Xi_y*F(f)):\Lambda_y\circ F(f)\Rightarrow G(f)\circ \Gamma_x,\]
for any morphism $f:x\to y$ in $\mc K$. When $F$ and $G$ are strict, this simplifies to 
\[G(f)* \Xi_x=\Xi_y* F(f).\] 

An \emph{equivalence} between the objects $x,y\in\mc K$ consists of two morphisms $f:x\leftrightarrow y:g$ together
with invertible 2-morphisms $\alpha:g\circ f\iso \id_x$ and $\beta:f\circ g\iso \id_y$.

If $F:\mc K\to \mc L$ is a 2-functor between 2-categories, then in general we say $F$
is ``locally $\ph$" if $\ph$ is a predicate applicable to functors between 1-categories 
which holds of each functor $\mc K(x,y)\to \mc L(F(x),F(y))$ induced by $F$. For instance, 
we can in this way ask that $F$ be locally essentially surjective, locally fully faithful, or 
locally an equivalence. 

We shall use the phrase ``bicategorically $\ph$" for global properties of $F$ that categorify
the property $\ph$ as applied to a single functor of categories. For instance, we shall use 
``bicategorically fully faithful" as a synonym for the potentially misleading term ``local equivalence". Finally, $F$ will be said to be ``bicategorically conservative" if 
it reflects equivalences: whenever we have $f:x\to y$ in $\mc K$ such that $F(f)$ is an 
equivalence in $\mc L$, we can conclude $f$ is an equivalence in $\mc K$.
\end{definition}

\begin{remark}\label{ffconservative}
Observe that any bicategorically fully faithful 2-functor 
$F:\mc K\to \mc L$ is bicategorically conservative. 
Given an equivalence $H(f):H(x)\leftrightarrows H(y):g$, 	
since quasi-inverses are closed under isomorphism and $H$ is locally essentially surjective,
we may assume that $g=H(g')$ for some $g':y\to x$. Now we have only to note that 
$g'f$ and $fg'$ are isomorphic to their respective identities since 
$H(g')H(f)$ and $H(f)H(g')$ are. 
\end{remark}

We now recall the definitions relevant to the theory of derivators.

\begin{definition}\label{HFin}
Denote by $\twocat{HFin}$ the 2-category of \emph{homotopy finite} categories. A category
is homotopy finite, often (confusingly) called finite direct, if its nerve has finitely many 
nondegenerate simplices; equivalently, if it is finite, skeletal, and admits no nontrivial 
endomorphisms.
\end{definition}

\begin{definition}
A \emph{prederivator} is a 2-functor $\D:\twocat{Dia}^{\op}\to\twocat{CAT}$
into the 2-category $\twocat{CAT}$ of large categories. The 2-category $\twocat{Dia}$ will be,
for us, either the 2-category of small categories $\twocat{Cat}$
or the 2-category $\twocat{HFin}$ of homotopy finite categories. 
\end{definition} 

We will often denote $\D(u)$ by $u^*$, for $u:J\to K$
a functor in~$\twocat{Dia}$, and similarly for a 2-morphism $\a$ in $\twocat{Dia}$. 

For categories $J,K\in \twocat{Dia}$, we have a functor $\dia_J^K:\D(J\times K)\to \D(J)^K$
induced by the action of $\D$ on the functors and natural
transformations from $[0]$ to $K$. We refer to $\dia_J^K$ as a ``partial underlying 
diagram functor," and when $J=[0]$ simply as the ``underlying diagram 
functor," denoted~$\dia^K$. See \cite[Section 1]{groth} for more on
the associated diagram.

Below are those axioms of derivators that are relevant to this paper. We
stick with the traditional numbering 
due to Maltsiniotis \cite{maltsiniotis} but leave out the axioms we 
shall not consider.
The 2-functor $\D$ is 
a \emph{semiderivator} if it satisfies the first two of the following
axioms, and \emph{strong} if it satisfies (Der5). We introduce here
a variant (Der5\textsf{\large'}) of the fifth axiom, prederivators satisfying
which will be called \emph{smothering}.
Let us remark that, in the presence of Axiom \hyperref[Der2]{(Der2)}, 
\hyperref[Der5']{(Der5\textsf{\large'})} requires exactly that $\dia_J^{[1]}$
be smothering in the sense of \cite{riehl}, which explains the nomenclature.

\begin{enumerate}[(Der1)]
\item\label{Der1}
Let $(J_i)_{i\in I}$ be a family of objects of $\twocat{Dia}$
such that $\coprod_I J_i\in\twocat{Dia}$.
Then the canonical map
$\D\left(\coprod_I J_i\right)\to\prod_I\D(J_i)$
is an equivalence.

\item\label{Der2}
For every $J\in\twocat{Dia}$, the underlying diagram functor
$\dia^J:\D(J)\to\D([0])^J$
is conservative, i.e., reflects equivalences.

%

\item[(Der5)]\label{Der5}
For every $J\in\twocat{Dia}$, the partial underlying
diagram functor $\dia_J^{[1]}:\D(J\times [1])\to\D(J)^{[1]}$
is full and essentially surjective on objects.

\item[(Der5\textsf{\large'})]\label{Der5'}
For every $J\in\twocat{Dia},$ the partial underlying
diagram functor $\dia_J^{[1]}:\D(J\times[1])\to \D(J)^{[1]}$ is full and surjective on objects.
\end{enumerate}

A morphism of prederivators is a pseudonatural transformation, and a 2-morphism
is a modification (see Definition \ref{2Cat}.) Altogether, we get the 2-category $\twocat{PDER}_{\twocat{Dia}}$ of prederivators defined on $\twocat{Dia}$. We
shall make use of the shorthand $\twocat{PDER}$ to represent a 2-category 
of prederivators defined on an arbitrary $\twocat{Dia}$.
 When we insist on strictly 2-natural transformations, we get the 
sub-2-category $\twocat{PDER}^\str$, of which we will primarily use the underlying
category, $\cat{PDER}^\str$. Our capitalization convention will reserve the notation 
$\twocat{PDer}$ for a 2-category of ``small" prederivators defined as those
 with values in the 2-category $\twocat{Cat}$ of small categories.


\section{The basic construction}

In this section, we will define the homotopy prederivator associated to a quasicategory
as the value of a functor, a 2-functor, and a simplicial functor. 

\subsection*{The prederivator associated to a quasicategory}
We denote the category associated to the poset $0<1<\cdots<n$ 
by $[n]$, so that $[0]$ is the terminal category.
The simplex category
$\Delta$ is the full subcategory of $\cat{Cat}$
on the categories $[n]$.

If $S$ is a simplicial set, that is, a 
functor $\Delta^{\op}\to\cat{Set}$, then we denote
its set of $n$-simplices by $S([n])=S_n$.
The face map $S_n\to S_{n-1}$ which forgets
the $i^{\mathrm{th}}$ vertex will be denoted
$d^n_i$ or just~$d_i$. We denote by $\Delta^n$ the simplicial set represented
by $[n]\in\Delta$. Equivalently, 
$\Delta^n=N([n])$, where we recall that the nerve $N(J)$ of a category~$J$ is the simplicial
set defined by the formula $N(J)_n=\cat{Cat}([n],J)$. The natural 
extension of $N$ to a functor is a fully faithful embedding of categories in 
simplicial sets. See \cite[Proposition~B.0.13]{joyal}.

We recall that a \emph{quasicategory} \cite{joyal},
 called an $\infty$-category in \cite{lurie}, is a simplicial set $Q$ in which every inner 
horn has a filler. That is, every map 
$\Lambda_i^n\to Q$ extends to an $n$-simplex 
$\Delta^n\to Q$ when $0<i<n$, where $\Lambda_i^n\subset \Delta^n$ is the simplicial subset 
generated by all faces $d_j\Delta^n$ with $j\neq i$. 
For instance, when $n=2$, the only inner horn is $\Lambda^2_1$, and then the 
filler condition simply says we may compose ``arrows" (that is, 1-simplices) 
in $Q$, though not uniquely. Morphisms of quasicategories are simply
morphisms of simplicial sets. The quasicategories
in which every inner horn has a \emph{unique} filler
are, up to isomorphism, the nerves of categories;
in particular the nerve functor $N:\cat{CAT}\to\cat{SSET}$ 
factors through the subcategory of quasicategories,~$\cat{QCAT}$.

\hypertarget{HoQ}{} Every quasicategory $Q$ has a homotopy category
$\Ho(Q)$, the 1-category defined as follows.
The objects of $\Ho(Q)$ are simply the 0-simplices of $Q$.
For two 0-simplices $q_1,q_2$, temporarily define
$Q_{q_1,q_2}\subset Q_1$ to be the set of 1-simplices
$f$ with initial vertex $q_1$ and final
vertex $q_2$. Then the hom-set
$\Ho(Q)(q_1,q_2)$ is the quotient of $Q_{q_1,q_2}$
which identifies \emph{homotopic} 1-simplices. Here 
two 1-simplices $f_1,f_2\in Q_{q_1,q_2}$ are said
to be homotopic if $f_1,f_2$ are two faces of some
2-simplex in which the third face is both outer and degenerate. 
We have a functor $\Ho:\cat{QCAT}\to\cat{CAT}$ 
from quasicategories to categories, left adjoint to the nerve 
$N:\cat{CAT}\to \cat{QCAT}$. This follows from the fact that 
a morphism $f:Q\to R$ of quasicategories preserves the homotopy 
relation between 1-simplices, so that it descends to a well defined functor 
$\Ho(f):\Ho(Q)\to \Ho(R)$. In fact, $\Ho:\cat{QCat}\to\cat{Cat}$ 
admits an extension, sometimes denoted $\tau_1$, to all of $\cat{SSet}$,
which is still left adjoint to $N$. But it is not amenable to computation.

The fact that the Joyal model structure is Cartesian and has the quasicategories 
as its the fibrant objects implies (see \cite[2.2.8]{riehl})
that $Q^S$ is a quasicategory for every
simplicial set $S$ and quasicategory $Q$.
In particular, the category of quasicategories is enriched over itself
via the usual simplicial exponential 
\[(R^Q)_n=\cat{SSET}(Q\times\Delta^n,R).\]

It is immediately checked that the homotopy category functor 
$\Ho$ preserves finite products, so that
by change of enrichment we get finally \emph{the 2-category of quasicategories}, 
$\twocat{QCAT}.$
Its objects are quasicategories,
and for quasicategories $Q,R,$ the 
hom-category~$\twocat{QCAT}(Q,R)$ is simply
the homotopy category $\Ho(R^Q)$ of the 
hom-quasicategory~$R^Q$. This permits the following tautological 
definition of equivalence of quasicategories.

\begin{definition}\label{equivalenceofquasicategories}
An equivalence of quasicategories is an equivalence
in~$\twocat{QCAT}$.
\end{definition}

\begin{remark}\label{2catequiv}
Thus an equivalence of quasicategories
is a pair of maps $f:Q\leftrightarrows R:g$ together with two homotopy
classes $a=[\a],b=[\b]$ of morphisms $\a:Q\to Q^{\Delta^1},
\b:R\to R^{\Delta^1}$, with endpoints $gf$ and $\id_Q$, respectively,
$fg$ and $\id_R$, such that $a$ is an isomorphism
in $\Ho(Q^Q)$, as is $b$ in $\Ho(R^R)$. We can make
the definition yet more explicit by noting that, for each
$q\in Q_0$, the map $\a$ sends $q$ to some $\a(q)\in Q_1$, 
and recalling that the invertibility of $a$ is equivalent to that
of each homotopy class $[\a(q)]$, as 
explicated for instance in the statement below:
\end{remark}

\begin{lemma}[\cite{riehl}, 2.3.10]\label{pointwiseness} 
The equivalence class $[\a]$ of a map $\a:Q\to R^{[1]}$ is an isomorphism
in the homotopy category $\Ho(R^Q)$ if and only if, for 
every vertex $q\in Q_0$ of~$Q$, the equivalence class
$[\a(q)]$ is an isomorphism in $\Ho(R)$.
\end{lemma}

We now construct the 2-functor $\HO:\twocat{QCAT}\to\twocat{PDer}$
(with respect to an arbitrary $\twocat{Dia}$.) Restricting to $\twocat{QCat}$
gives us all the forms of $\HO$ of interest to us. 

We first extend $\Ho$ to a 2-functor of the same name,
$\Ho:\twocat{QCAT}\to\twocat{CAT}$. 
This still sends a quasicategory to its homotopy 
category; we must define the action
on morphism categories. This will be for each $R$ and $Q$ a functor
\[\Ho_{Q,R}:\twocat{QCAT}(Q,R)=\Ho(R^Q)\to
\Ho(R)^{\Ho(Q)}=\twocat{CAT}(\Ho(Q),\Ho(R)).\]

The functor  $\Ho_{Q,R}$ is defined as the
transpose of the composition 

\[\Ho(R^Q)\times \Ho(Q)\iso\Ho(R^Q\times Q)\map{\Ho(\mathrm{ev})}\Ho(R)\]
across the product-hom adjunction in the
1-category~$\cat{Cat}$.
For this isomorphism we have used again the preservation of finite products
by $\Ho$. The morphism $\mathrm{ev}:R^Q\times Q\to R$
is evaluation, the counit of the adjunction
$(-)\times Q\ladj (-)^Q$ between endofunctors of~$\cat{QCAT}$.  

We also need a 2-functor
 $N:\twocat{CAT}\to\twocat{QCAT}$ 
sending a category $J\in\twocat{CAT}$
to $N(J)$. The map on hom-categories is
the composition
$J^K\iso \Ho(N(J^K)) \iso \Ho(N(J)^{N(K)})$.
The first isomorphism is the inverse of the counit
of the adjunction $\Ho\ladj N$, which is an 
isomorphism by
full faithfulness of the nerve.
The second uses the fact that $N$ preserves exponentials, see 
\cite[Proposition~B.0.16]{joyal}.

Finally, we require the following 
fact: a monoidal functor $F:\mc V\to \mc W$
induces a 2-functor $F_*(-):\mc V-\twocat{Cat}\to\mc W-\twocat{Cat}$ between
2-categories of $\mc V$- and $\mc W$-enriched 
categories. 
The fully general version of this claim was apparently not published until recently; 
it comprises Chapter 4 of \cite{cruttwell}.
In our case, the functor $\Ho$ is monoidal insofar as it preserves products 
and thus it induces the 2-functor $\Ho_*(-)$ sending simplicially
enriched categories, simplicial functors, and simplicial natural transformations to
2-categories, 2-functors, and 2-natural transformations.

Now we define the homotopy prederivator. 

\begin{definition}\label{defofHOQ}
Let $Q$ be a quasicategory. Then the \emph{homotopy prederivator}
$\HO(Q):\twocat{Dia}^{\op}\to\twocat{CAT}$ is given as
the composition
\[\twocat{Dia}^{\op}\map{N^{\mathrm{op}}}
\twocat{QCAT}^{\op}
\map{Q^{(-)}}\twocat{QCAT}
\map{\Ho}\twocat{CAT}.\]

\end{definition}

In particular, $\HO(Q)$ maps a category~$J$ to 
the homotopy category of $J$-shaped diagrams in
$Q$, that is, to~$\Ho(Q^{N(J)})$. Given a morphism of quasicategories $f:Q\to R$,
we have a strictly 2-natural morphism of prederivators (see Definition \ref{2Cat})
$\HO(f):\HO(Q)\to\HO(R)$ given as the analogous
composition $\HO(f)=\Ho\circ f^{(-)}\circ N$,
so that for each category $J$ the functor $\HO(f)_J$ is given by
post-composition with $f$, that is, by
$\Ho(f^{N(J)}):\Ho(Q^{N(J)})\to\Ho(R^{N(J)})$.

We now record the axioms which are satisfied by the homotopy prederivator of any quasicategory. 
First, a quasicategorical lemma:

\begin{lemma}\label{liftingsquares}
Let $Q$ be a quasicategory, and $X:[1]\times[1]\to\Ho(Q)$ a commutative 
square in its homotopy category. Suppose we have chosen $f,g\in Q_1$ representing
the vertical edges of $X$, so that $[f]=X|_{\{0\}\times[1]}$ and 
$[g]=X|_{\{1\}\times[1]}$. Then there exists $\widehat X:f\to g$ in 
$\Ho(Q^{\Delta^1})$ lifting $X$, in the sense that 
$0^*\widehat X=X|_{[1]\times\{0\}}$ and $1^*\widehat X=X|_{[1]\times\{1\}}$.
\end{lemma}

\begin{proof}
We must show that any homotopy-commutative
square $X:[1]\times[1]\to\Ho(Q)$ with chosen lifts $f,g\in (Q)_1$ of its left and right edges 
underlies a morphism $\widehat X:f\to g$ in $\Ho(Q^{\Delta^1})$. 
For this we first lift the top and bottom edges of $X$ to some 
$h$ and $k$ in $Q_1$ and choose 2-simplices $a,b$ filling the horns
\[\begin{tikzcd}
\bullet\ar[r,"h"]&\bullet\ar[d,"g"]&\bullet\ar[d,"f"]&\\
&\bullet &\bullet\ar[r,"k"]&\bullet 
\end{tikzcd}\]
respectively. Let $m=d_1 a$ and $n=d_1 b$ be the inner faces of $a$ and $b$.

Since $X$ was homotopy commutative, we know $[g]\circ[h]=[k]\circ[f]$ in $\Ho(Q)$, 
that is, $[m]=[n]$.
So there exists a 2-simplex $c$ with boundary 
\[\begin{tikzcd}
\bullet\ar[dr,"m"]\ar[rr,equals]&&\bullet\ar[dl,"n"]\\
&\bullet
\end{tikzcd}\]
giving a homotopy
between $d_1a$ and $d_1b$. 

Now we define a map $H:\Lambda^3_1\to Q$ with  
$d_0H=b,d_2H=c$, and $d_3H$ degenerate on $f$, as below:
\[\begin{tikzcd}
&\bullet\ar[dr,"f"]\ar[rr,"n"]&&\bullet &&\bullet\ar[rr,"n"]&&\bullet \\
\bullet\ar[rr,"f"]\ar[ur,equals]&&\bullet\ar[ur,"k"]&&\bullet\ar[rr,"f"]\ar[ur,equals]\ar[urrr,"m"]&&\bullet\ar[ur,"k"]
\end{tikzcd}\]

Filling $H$ to a 3-simplex $\hat H$, we get the desired square by juxtaposing 
$d_1\hat H$ and $a$. 
\end{proof}

\begin{prop}
For any quasicategory $Q$, the homotopy prederivator 
$\HO(Q)$ satisfies the axioms \hyperref[Der1]{(Der1)},
\hyperref[Der2]{(Der2)}, and \hyperref[Der5']{(Der5\textsf{\large'})}.
\end{prop}

\begin{proof}
The axiom \hyperref[Der1]{(Der1)} follows from the fact that $Q\mapsto Q^J$ preserves coproducts in
$J$, and that $\Ho$ preserves all products. \hyperref[Der2]{(Der2)} is an application of 
Lemma \ref{pointwiseness}, with $Q$ specialized to $N(J)$ for some $J\in\twocat{Dia}$.
For \hyperref[Der5']{(Der5\textsf{\large'})}, surjectivity of $\dia_J^{[1]}$ follows immediately
from the definition of the homotopy category. Fullness is exactly Lemma \ref{liftingsquares}.
\end{proof}

It may be worth noting that, while it is possible to define a 2-category 
$\twocat{SSet}$ of simplicial sets using $\tau_1$ and extend 
$\HO$ to $\twocat{SSet}$, the prederivator associated to an arbitrary simplicial 
set $S$ will not, in general, satisfy any of the three axioms. 
It is straightforward to see that $\HO(S)$ need not satisfy 
\hyperref[Der2]{(Der2)} or \hyperref[Der5]{(Der5)}, while the reason 
\hyperref[Der1]{(Der1)} may fail is that $\tau_1$, unlike $\Ho$, need
not preserve infinite products.

\subsection*{The simplicial enrichment of prederivators}

The 2-functor $\HO:\twocat{QCAT}\to \twocat{PDer}$ factors
through the subcategory $\twocat{PDer}^\str$ 
in which the morphisms are required to be strictly
2-natural. Its underlying category $\cat{PDer}^\str$ admits 
a simplicial enrichment $\scat{PDer}$, as we now recall.

Muro and Raptis showed how to define the simplicially enriched category 
$\scat{PDer}$ in \cite{muro}. First, note that for
 any prederivator $\D$ and each category $J\in \twocat{Dia}$ we have a shifted 
prederivator $\D^J=\D\circ(J\times -)$. 
This shift is a special case of the cartesian closed structure on 
$\twocat{PDer}$ discussed in \cite[Section 4]{heller2}.
Explicitly, given two prederivators $\D_1,\D_2,$ and denoting by
$\widehat J$ the prederivator represented by a small category~$J$,
the exponential is defined by $\D_2^{\D_1}(J)=\twocat{PDer}(\widehat J \times \D_1,\D_2)$.
Then the 2-categorical Yoneda lemma implies that 
the shifted prederivator $\D^J$ is canonically isomorphic to  the 
prederivator exponential $\D^{\widehat J}$. This allows us to interpret
expressions such as $\D^\a:\D^u\Rightarrow\D^v:\D^K\to\D^J$, 
when $\a:u\Rightarrow v:J\to K$ is a natural transformation,
by using the internal hom 2-functor.

\begin{remark}\label{simp}
For a natural transformation $\a:u\Rightarrow v:J\to K$
between functors in $\cat{Cat}$, the preceding definition of $\D^\a$ gives only
a shadow of the full action of $\a$ on $\D$.
The natural transformation $\a$ corresponds naturally to a functor
$\bar\a: J\times [1]\to K$, associated to which
we have a prederivator morphism
$\D^{\bar\a}:\D^K\to \D^{J\times [1]}$,
that is, a family of functors
$\D(K\times I)\to\D(J\times I\times [1])$.
This is strictly more information, as composing
with 
the underlying diagram functor
$\dia^{[1]}_{J\times I}:\D(J\times I\times [1])\to \D(J\times I)^{[1]}$
recovers our original $\D^\a$.
What is happening here is that the entity
$\D^{(-)}$ is more than a 2-functor 
$\twocat{Cat}^{\op}\to\twocat{PDer}$:
it is a simplicial functor 
$N_*\twocat{Cat}^{\op}\to\scat{PDer}$ from the 
simplicial category of nerves of categories
to the simplicial category of prederivators,
which we must now define.
\end{remark}

For each category $J$ let $\diag_J:J\to J\times J$
be the diagonal functor. 
\begin{definition}
We define 
$\scat{PDer}$ as a simplicially enriched category 
whose objects are the prederivators. 
The mapping simplicial 
sets have $n-$simplices as follows:
$\cat{PDer}_n(\D_1,\D_2)=\cat{PDer}^\str(\D_1,\D^{[n]}_2)$.
For $(f,g)\in \cat{PDer}_n(\D_2,\D_3)\times\cat{PDer}_n(\D_1,\D_2)$, the composition $f * g:\D_1\to \D_3^{[n]}$ is given by the formula below,
in which we repeatedly apply the internal hom
2-functor discussed above Remark \ref{simp}.
 \[\label{MRcomposition}\D_1\map{g} \D_2^{[n]}\map{f^{[n]}}
 \left(\D_3^{[n]}\right)^{[n]}\iso\D_3^{[n]\times [n]}\map{\D_3^{\diag_{[n]}}}\D_3^{[n]}\]
\end{definition}

We can now extend $\HO$ to a simplicially enriched functor.
The definition follows formally from the following 
interpretation of the simplicial enrichments on
$\cat{QCAT}$ and $\cat{PDer}^\str$. Each category
has a given cosimplicial object, respectively given by the representable simplicial sets
$\Delta^\bullet$ and the representable prederivators $\widehat{[\bullet]}$. We have a 
natural isomorphism
$\widehat{[\bullet]}\iso\HO(\Delta^\bullet)$ following from the full faithfulness of the nerve. 
This shows that for any quasicategory $R$, the simplicial prederivator $\HO(R)^{[\bullet]}$ is 
isomorphic to $\HO(R^{\Delta^\bullet})$. In particular, we can define 
$\HO:\scat{QCAT}(Q,R)\to \scat{PDer}(\HO(Q),\HO(R))$ on $n$-simplices by the composition
\[\cat{QCAT}(Q,R^{\Delta^n})\to\cat{PDer}^\str(\HO(Q),\HO(R^{\Delta^n}))\iso \cat{PDer}^\str(\HO(Q),\HO(R)^{[n]}).\]
We also have canonical cosimplicial objects in $\cat{QCAT}$ and $\cat{PDer}^\str$, 
given by $\Delta^\bullet\times \Delta^\bullet$ and $\widehat{[\bullet]}\times \widehat{[\bullet]}$ 
respectively, which are mapped isomorphically into each other by $\HO$. 
Furthermore, the isomorphisms
$\HO(\Delta^\bullet)\iso \widehat{[\bullet]}$ and $\HO(\Delta^\bullet\times\Delta^\bullet)\iso 
\widehat{[\bullet]}\times\widehat{[\bullet]}$ commute with the diagonal and projection morphisms
used in the definition of the simplicial compositions in $\scat{QCAT}$ and $\scat{PDer}$, 
so that the suggested definition $\HO$ respects the simplicial compositions.

In \cite{muro} a restriction of this 
enrichment, which we now recall, was of primary interest.
Each prederivator $\D$ has an ``essentially constant"
shift by a small category $J$ denoted $\D_{\eq}^J$.
This is defined as follows:
$\D_\eq^J(K)\subset \D(J\times K)$ is the full
subcategory on those objects $X\in \D(J\times K)$
such that in the partial underlying diagram 
$\dia_K^J (X)\in \D(K)^J$, the image of every morphism
of $J$ is an isomorphism in $\D(K)$. We shall only
need $J=[n]$, when an object of $\D_\eq^{[n]}(K)$
has as its partial underlying diagram
 a chain of $n$ isomorphisms in~$\D(K)$.

Then we get another simplicial enrichment:

\begin{definition}\label{scateq}
The simplicial category $\scat{PDer}^\eq$
is the sub-simplicial category of $\scat{PDer}$
with 
\[\cat{PDer}^\eq_n(\D_1,\D_2)=\cat{PDer}^\str(\D_1,\D_{2,\eq}^{[n]}).\]
\end{definition}
This leads to the notion of equivalence of prederivators under which Muro and 
Raptis showed Waldhausen K-theory is invariant. We say \emph{coherent}
below where Muro and Raptis use \emph{strong}, to avoid ambiguity.

\begin{definition}\label{coherentequivalence}
A coherent isomorphism between strict 2-natural transformations
$F:\D\to \mathscr{E},G:\mathscr{E}\to\D$ between prederivators is given by 
2-morphisms $\a:\D\to\D_\eq^{[1]}$ and \mbox{$\b:\mathscr{E}\to\mathscr{E}_\eq^{[1]}$} 
such that the vertices of $\a$ are $GF$ and $\id_{\D}$, and similarly for 
$\b$. 

The 2-natural transformations $F$ and $G$ comprise a \emph{coherent equivalence}
of prederivators if there exist zigzags of coherent isomorphisms connecting
$GF$ to $\id_{\D_1}$ and $FG$ to $\id_{\D_2}$. 
\end{definition}

By \cite[Proposition~B.0.15]{joyal},
the extension $\tau_1:\cat{SSet}\to\cat{Cat}$ of the homotopy category
functor $\Ho$ to the entirety of $\cat{SSet}$ preserves finite products.
Thus the simplicial categories $\scat{PDer}$ and $\scat{PDer}^\eq$ 
give rise to 2-categories $\tau_{1*}\twocat{PDer}$ and 
$\tau_{1*}\twocat{PDer}^\eq$ by applying $\tau_1$ to each
hom-simplicial set.

\begin{remark}\label{rem:typesofEquivs}
We are now provided with an abundance of notions of equivalences of
prederivators. To wit, we have:
\begin{enumerate}[(1)]
\item The equivalences in the 2-category $\tau_{1*}\scat{PDer}^\eq$.
\item The equivalences in the 2-category $\tau_{1*}\scat{PDer}$.
\item The coherent equivalences as defined above.
\item The equivalences in the 2-category $\twocat{PDer}^\str$.
\item The morphisms in $\twocat{PDer}^\str$ which induce levelwise 
equivalences of categories.
\item The equivalences in the 2-category $\twocat{PDer}$. 
\item The morphisms in $\twocat{PDer}$ which induce levelwise 
equivalences of categories.
\end{enumerate}
\end{remark}

\begin{prop}\label{prop:typesofEquivs}
The following implications hold among the above classes of equivalences.
\[\begin{tikzcd}
(1)\ar[r,Rightarrow]\ar[d,Rightarrow]&(2)\ar[d,Rightarrow]\\
(3)\ar[r,Rightarrow]&(4)\ar[r,Rightarrow]&(5)\ar[r,Rightarrow]&(6)\ar[r,Leftrightarrow]&(7)
\end{tikzcd}\]

\end{prop}
\begin{proof}
The implications $(4)\implies (5)$ and $(6)\implies (7)$ are immediate,
since evaluation at any $J$ gives a 2-functor $\twocat{PDer}\to \twocat{Cat}$.
Since $(5)\implies (7)$ trivially, to prove $(5)\implies (6)$ it is enough to
show $(7)\implies (6)$. Given a pseudonatural transformation $F:\D_1\to \D_2$
inducing equivalences of categories $F_J:\D_1(J)\to \D_2(J)$ for all 
$J\in\twocat{Dia}$, arbitrary choices of quasi-inverses $G_J:\D_2(J)\to \D_1(J)$
to the $F_J$ can be compiled into a quasi-inverse for $F$ in $\twocat{PDer}$,
as can easily be checked directly.
More abstractly, this is a case of a basic result in 2-dimensional universal 
algebra, since $\twocat{PDer}^\str$ is the 2-category of strict algebras for a 
2-monad on the 2-category of $\twocat{Dia}$-indexed families of categories. See
\cite[Proposition~4.10]{lack}.

It thus remains to handle the implications involving the simplicial 
enrichments. Recall that if $\mc K_\bullet$ is a simplicial category, then the 
$2$-category $\tau_{1*}\mc K$ has as hom-categories $\mc K_{\tau_1}(x,y)=\tau_1(\mc K(x,y).$ 
Thus an equivalence in $\tau_{1*}\mc K$ is given by morphisms 
$f:x\leftrightarrows y:g$ together with four strings of edges in hom-simplicial 
sets of $\mc K_\bullet$, linking $\id_x$ and $\id_y$ to and from $gf$ and $fg$, 
respectively, which become mutually inverse in $\tau_{1*}\mc K(x,x)$, 
respectively $\tau_{1*}\mc K(y,y)$. 

This description yields the implications $(1)\implies(2)$, since an edge in
$\scat{PDer}^\eq$ is an edge in $\scat{PDer}$, and 
$(1)\implies (3)$, since a string of edges in $\scat{PDer}^\eq$ is in 
particular a zigzag. For $(3)\implies (4)$, consider a coherent equivalence
$F:\D_1\leftrightarrows \D_2:G$. Then by assumption, we may choose 
a zigzag \[\id_{\D_1}\to T_1 \leftarrow T_2\to...\to T_n\leftarrow GF\]
in which the arrows are 1-simplices in $\scat{PDer}^\eq(\D_1,\D_1)$. By
definition of $\scat{PDer}^\eq$, each such 1-simplex induces an invertible
modification in $\twocat{PDer}^\str(\D_1,\D_1)$. Thus we may compose along
the zigzag, inverting as necessary, to get a modification $\id_{\D_1}\iso GF$. 
Similarly, we get $FG\iso \id_{\D_2}$, so $F$ and $G$ are mutually quasi-inverse
in $\twocat{PDer}^\str$. 
For $(2)\implies (4)$, the natural morphism of simplicial sets 
$\scat{PDer}(\D_1,\D_1)\to N\twocat{PDer}^\str(\D_1,\D_1)$
induces an identity-on-objects functor 
$\tau_{1*}\scat{PDer}(\D_1,\D_1)\to \twocat{PDer}^\str(\D_1,\D_1)$. 
Thus if $\id_{\D_1}$ and $GF$ are isomorphic 
in $\tau_{1*}\scat{PDer}(\D_1,\D_1)$, they are isomorphic in $
\twocat{PDer}^\str(\D_1,\D_1)$. 

\end{proof}

\begin{remark}\label{modificationsarebad}
In general, no implications between the various classes of equivalences
hold other than those listed above. We forebear to give explicit counterexamples, 
but suggest the nature of each obstruction here. To get $(6)\implies (5)$ or 
$(5)\implies (4)$ would require a strictifiability theorem for pseudonatural 
transformations, which does not generally hold. 
The implications $(4)\implies (3)$ and $(4)\implies (2)$ are versions of 
Axiom \ref{Der5}, weakened to allow lifting into strings or zigzags of coherent 
morphisms but strengthened to apply to endomorphisms of $\D$, rather than just to 
objects in the values of $\D$. 
The easiest way to get $(3)\implies (1)$ is for $\scat{PDer}(\D_i,\D_i)$ to admit
composition of edges, while the easiest way to get $(2)\implies (1)$ is for 
$\scat{PDer}^\eq(\D_i,\D_i)$ to admit inversion of 1-arrows. 

As the above remarks suggest, many of these implications collapse when the 
$\D_i$ are associated to quasicategories. 
By Theorem \ref{Maintheorem} below, any \emph{strict}
2-natural transformation $F:\HO(Q)\to\HO(R)$ arises from a morphism $f:Q\to R$ of 
quasicategories. By Theorem \ref{whiteheadforquasicats}, a map 
$f:Q\to R$ of quasicategories is
an equivalence if and only if $\HO(f)$ is an equivalence in $\twocat{PDer}$. 
Furthermore, the simplicial clause of Theorem \ref{Maintheorem} shows that an 
equivalence of quasicategories induces an equivalence in 
$\tau_{1*}\scat{PDer}^\eq$. Thus, for prederivators associated
to quasicategories, all of $(1),(2),(3),(4),$ and $(5)$ are equivalent. 

When the
quasicategories $Q$ and $R$ are small and $\twocat{Dia}=\twocat{Cat}$, Theorem 
\ref{Mainthm2cat} implies that even if $F:\HO(Q)\to\HO(R)$ is merely pseudonatural, if
it is an equivalence in $\twocat{PDer}$ then it arises from an equivalence of quasicategories.
Thus in that case, all seven notions of equivalences are equivalent. However, in Conjecture 
\ref{conj:incoherent} we suggest a pseudonatural equivalence between $\twocat{HFin}$-indexed 
prederivators of small quasicategories which does not arise from an equivalence, so we
do not expect a full collapse of equivalence notions in general, even for quasicategories.
\end{remark}


\section{The simplicial embedding $\scat{QCAT}\to\scat{PDer}$}\label{section:simplicial}
In this section, we prove that categories of 
arbitrarily large quasicategories embed fully faithfully in any category
of prederivators and strict morphisms. We extend this result to a fully faithful
embedding of simplicial categories, as well as of categories enriched in Kan 
complexes.

\begin{thm}\label{Maintheorem} 
The ordinary functor $\HO:\cat{QCAT}\to\cat{PDer}^\str$ is fully faithful.
\end{thm}

We first give some corollaries.

\begin{cor}\label{simplicialembedding}
The simplicial functor $\HO:\scat{QCAT}\to\scat{PDer}$
is simplicially fully faithful.	
\end{cor}
\begin{proof}
The action of $\HO$ on $n$-simplices was defined as map of
$\cat{QCAT}(Q,R^{\Delta^n})\to \cat{PDer}^\str(\HO(Q),\HO(R^{\Delta^n}))$ induced
by $\HO$ 
followed by the canonical isomorphism 
\[\cat{PDer}^\str(\HO(Q),\HO(R^{\Delta^n}))\iso \cat{PDer}^\str(\HO(Q),\HO(R)^{[n]}).\] 
Thus, $\HO$ induces an isomorphism on $n$-simplices for every $n$ and is simplicially fully 
faithful.
\end{proof}

Define, for the moment,
$\scat{QPDer}\subset\scat{PDer}$ to be the image of quasicategories
in prederivators, so that Corollary \ref{simplicialembedding} gives an isomorphism of 
simplicial categories $\scat{QCAT}\iso \scat{QPDer}$. In particular,
$\scat{QPDer}$ is not merely a simplicial category, but actually a category
enriched in quasicategories.

Recall that the inclusion of Kan complexes into quasicategories
 has a right adjoint $\iota$, which we will call
the Kan core. For a quasicategory $Q$, the core 
$\iota Q$ is the sub-simplicial set such that an $n$-simplex
$x\in Q_n$ is in $(\iota Q)_n$ if and only if every 
1-simplex of $x$ is an isomorphism in $\Ho(Q)$. See~\cite[Section 1]{joyal2}.

As a right adjoint, $\iota$ preserves products,
so that for any quasicategorically enriched category $\mc C$ 
we have an associated Kan complex-enriched category~
$\iota_*\mc C$, given by taking the core homwise. (This change of enrichment
does not exist on a point-set level for general simplicially enriched categories,
which explains our inelegant introduction of $\scat{QPDer}$.)

\begin{cor}\label{bigdiagramcorollary}
The associated prederivator functor $\HO:\scat{QCAT}\to\scat{PDer}$
induces an isomorphism of $\cat{Kan}$-enriched categories
$\iota_*\HO:\iota_*\mathbf{QCAT}_{\bullet}\to \iota_*\mathbf{QPDer}_{\bullet}$.
\end{cor}
\begin{proof}
The given $\cat{Kan}$-enriched functor exists by the argument of \cite{cruttwell} described
above. It is defined predictably, in the manner
of Equation \ref{kanff} below.
We just have to show that $\iota_*\HO$
induces isomorphisms on hom-objects, since $\iota_*\HO$ is bijective on
objects by definition. 
Given the isomorphism 
$\HO_{Q,R}:\scat{QCAT}(Q,R)\iso\scat{PDer}(\HO(Q),\HO(R))$
of Theorem \ref{Maintheorem}, we get isomorphisms 
\begin{equation}\label{kanff}\iota(\HO_{Q,R}):\iota(\scat{QCAT}(Q,R))\iso\iota(\scat{PDer}(\HO(Q),\HO(R)))\end{equation}
as desired. 
\end{proof}

\begin{remark}
The $\cat{Kan}$-enriched category
$\iota_*\mathbf{QCAT}_{\bullet}$ is a model of the homotopy
theory of homotopy theories, which thus embeds
into prederivators. In particular, the homotopy category of homotopy theories
embeds in the simplicial homotopy category of $\scat{PDer}^\eq$.

In Section \ref{section:2cat}, we improve this to show that
\emph{the homotopy 2-category}
in the sense of \cite{riehl2} embeds in the 2-category $\twocat{PDer}_{\twocat{Cat}}$,
a much more concrete object, under certain size assumptions. 
The word \emph{the} is partially justified here by work of
Low \cite{Zhen} indicating that the 2-category $\twocat{QCat}$ has a universal role analogous to that of 
``\emph{the} homotopy category", namely, the homotopy category of spaces.
\end{remark}

We turn to the proof of Theorem \ref{Maintheorem}. 
We must show that the ordinary functor $\HO$ gives 
an isomorphism between the sets
$\cat{QCAT}(Q,R)$ and $\cat{PDer}^{\mathrm{str}}(\HO(Q),\HO(R))$.
This is Proposition~\ref{main prop}, whose proof has the following outline:
\begin{enumerate}[(1)]
\item\label{proofstep1} Eliminate most of the data of a
prederivator map by showing strict 
maps $\HO(Q)\to \HO(R)$ are determined
by their restriction to natural transformations
between ordinary functors 
$\cat{Cat}^{\op}\to\cat{Set}$. This is Lemma~\ref{catset lemma}.

\item\label{proofstep2} Show that 
$\HO(Q)$ and $\HO(R)$ recover
$Q$ and $R$ upon restricting the domain to 
$\Delta^{\op}$ and the codomain to $\cat{Set}$,
and that natural transformations as in the previous
step are in bijection with maps $Q\to R$. This
is Lemma~\ref{kan lemma}.

\item\label{proofstep3} Show that $\HO(f)$ restricts back to
$f$ for a map $f:Q\to R$, which implies that
$\HO$ is faithful, and that a map $F:\HO(Q)\to\HO(R)$
is exactly $\HO$ applied to its restriction, which implies
that $\HO$ is full. This constitutes the proof
of Proposition~\ref{main prop} proper.
\end{enumerate}

Let us begin with step \hyperref[proofstep1]{(1)}.
\begin{definition}\label{catset def}
A $\cat{Dia}$-set is a large presheaf on $\cat{Dia}$ 
that is, an ordinary functor
$\cat{Dia}^{\op}\to\cat{SET}$.
\end{definition}

Given a prederivator $\D$, let 
$\D^{\ob}:\cat{Dia}^{\op}\to\cat{SET}$ 
be its underlying $\cat{Dia}$-set, so that 
$\D^{\ob}$
sends a small category~$J$ to the set of objects
$\ob(\D(J))$ and a functor $u:I\to J$ to
the action of $\D(u)$ on objects.

Recall that where \hyperref[Der5]{(Der5)} requires that
$\dia:\D(J\times[1])\to\D(J)^{[1]}$ be
(full and) essentially surjective, \hyperref[Der5']{(Der5\textsf{\large'})}
insists on actual surjectivity on objects.
The following lemma shows that under this assumption
most of the apparent structure of a strict 
prederivator map is redundant.

\begin{lemma}\label{catset lemma}
A strict morphism $F:\D_1\to\D_2$ between
prederivators satisfying \hyperref[Der5']{(Der5\textsf{\large'})}
is determined by its restriction to the
underlying $\cat{Dia}$-sets 
$\D_1^{\ob},\D_2^{\ob}$. That is, the restriction 
functor from prederivators satisfying \hyperref[Der5']{(Der5\textsf{\large'})} to $\cat{Dia}$-sets
is faithful.
\end{lemma}

\begin{proof}
The data of a strict morphism $F:\D_1\to\D_2$
is that of a functor
$F_J:\D_1(J)\to\D_2(J)$ for every $J$.
\footnote{Note the simplification here over 
pseudonatural transformations, which require also a 
natural transformation associated to 
every functor and do not induce maps of
$\cat{Dia}$-sets. That is the fundamental difficulty leading to the dramatically 
different techniques of the next sections.}

The induced map $F^{\ob}:\D_1^{\ob}\to\D_2^{\ob}$
is given by the action of $F$ on objects. 
So to show faithfulness it is enough to show that,
given a family of functions
$r_J:\ob(\D_1(J))\to\ob(\D_2(J))$, that is,
the data required in a natural transformation between
$\cat{Dia}$-sets,
there is at most one 2-natural transformation
with components
$F_J:\D_1(J)\to\D_2(J)$ and object parts 
$\ob(F_J)=r_J$.

Indeed, suppose $F$ is given with object parts
$r_J=\ob(F_J)$ and
let $f:X\to Y$ be a morphism in
$\D_1(J)$. Then by Axiom \hyperref[Der5']{(Der5\textsf{\large'})}, $f$ is the
underlying diagram of some 
$\widehat f~\in~\D_1(J\times[1])$.  
By 2-naturality, the following square must
commute:
\[\begin{CD}
\D_1(J\times[1])@>F_{J\times[1]}>>\D_2(J\times [1])\\
@V\dia_J^{[1]} VV @V\dia_J^{[1]} VV\\
\D_1(J)^{[1]}@>F_J>> \D_2(J)^{[1]}
\end{CD}
\]
Indeed, $\dia_J^{[1]}$ is the action of a prederivator
on the unique natural transformation between
the two functors $0,1:[0]\to[1]$ from the
terminal category to the arrow category, as is described
in full detail below \cite[Proposition~1.7]{groth}.
Thus the square above is an instance of the axiom of 
respect for 2-morphisms. It follows that we must have
$F_J(f)=F_J(\dia_J^{[1]}\widehat f)=\dia_J^{[1]}(r_{J
\times [1]}(\widehat f))$.

Thus if $F,G$ are two strict morphisms $\D_1\to\D_2$
with the same restrictions to the underlying
$\cat{Dia}$-sets, they must coincide, as claimed.
\end{proof}

Note the above does not claim that the restriction
functor is full: the structure of a strict prederivator
map is determined by the action on objects of each $\D_1(J),\D_2(J),$ 
but it is not generally true that an arbitrary map of $\cat{Dia}-$sets 
will admit a well defined extension to morphisms.

We proceed to step \hyperref[proofstep2]{(2)} of the proof. 

Let us recall
the theory of pointwise 
Kan extensions for 1-categories. Let 
$F:\mc C\to\mc D$ and $G:\mc C\to \mc E$ be functors.
At least if $\mc C$ and $\mc D$ are small and $\mc E$ is complete, then we always
have a right Kan extension $F_*G:\mc D\to \mc E$
characterized by the adjunction formula
$\mc E^{\mc D}(H,F_*G)\iso 
\mc E^{\mc C}(H\circ F,G)$ and computed on objects
by 
\begin{equation}F_*G(d)=\lim_{d\downarrow F}G\circ q\end{equation}
Here
$d\downarrow F$ is the comma category with objects $(c,f:d\to F(c))$
and morphisms the maps in $\mc C$ making the appropriate triangle commute,
and $q:d\downarrow F\to \mc C$ is the projection.

\begin{lemma}\label{kan lemma}
Let $j:\Delta^{\op}\to\cat{Dia}^{\op}$ be the
inclusion. Then for any quasicategory~$R$, the $\cat{Dia}$-set
$\HO(R)^{\ob}$ underlying $\HO(R)$
is the right Kan extension of $R$
along~$j$.
\end{lemma}

\begin{proof}
For any small category $J$,
the $\cat{Dia}$-set $\HO(R)^{\ob}$ takes $J$
to the set of simplicial set maps from $J$ to $R$:
\[\HO(R)^{\ob}(J)=\ob(\Ho(R^{N(J)}))
=\cat{SSET}(N(J),R)\]
We shall show that the latter is the value required of
$j_*R$ at $J$, which exists and is calculated via Equation 5.1 since
$\cat{SET}$ is complete (in the sense of a universe in which its objects constitute
the small sets.)
 
First, one of the basic properties of presheaf
categories
implies that $N(J)$ is a colimit over its category
of simplices. That is, 
\mbox{$N(J)=\colim\limits_{\Delta\downarrow NJ} y\circ q$}, where
$q:\Delta\downarrow NJ\to \Delta$ is the projection
and $y:\Delta\to\cat{SSet}$ is the Yoneda embedding.

Then we can rewrite the values of 
$\HO(R)^{\ob}$ as follows:
\[\HO(R)^{\ob}(J)=\cat{SSET}(N(J),R)=\cat{SSET}(\colim\limits_{\Delta\downarrow NJ} y\circ q,R)\iso\]\[\lim\limits_{(\Delta\downarrow NJ)^{\op}}\cat{SSET}(y\circ q,R)\iso
\lim_{(\Delta\downarrow NJ)^{\op}}R\circ q^{\op}\]
The last isomorphism follows from the Yoneda
lemma.

The indexing category $(\Delta\downarrow N(J))^{\op}$
has as objects pairs $(n,f:\Delta^n\to N(J))$
and as morphisms $\bar a:(n,f)\to(m,g),$
the maps $a:\Delta^m\to \Delta^n$ such that
$f\circ a=g$. That is,
$(\Delta\downarrow N(J))^{\op}\iso N(J)\downarrow\Delta^{\op}$,
where on the right-hand side $N(J)$ is viewed
as an object of $\cat{SSET}^{\op}$.
Using the full faithfulness
of the nerve functor $N$, we see 
$(\Delta\downarrow N(J))^{\op}\iso J\downarrow\Delta^{\op}$,
where again $J\in \cat{Dia}^{\op}$.

Thus, if $q^{\op}$ serves also to name the
projection $J/\Delta^{\op}\to \Delta^{\op}$,
we may continue the computation above with
\[\HO(R)^{\ob}(N(J))\iso \lim\limits_{N(J)\downarrow\Delta^{\op}}R\circ q^{\op}\]
This is exactly the formula for $j_*R(J)$ recalled
above. The isomorphism thus constructed is certainly
natural with respect to the action on maps of the Kan extension, 
so the lemma is established.
\end{proof}

We arrive at step \hyperref[proofstep3]{(3)}.

\begin{prop}\label{main prop}
The homotopy prederivator functor $\HO:\cat{QCAT}\to\cat{PDer}^\str$ is a fully faithful embedding 
of 1-categories.
\end{prop}

\begin{proof}

Note that, by Lemma~\ref{kan lemma}, 
the restriction of $\HO(Q)^\ob$ to a
functor $\Delta^{\op}\to\cat{SET}$
is canonically isomorphic to $Q$, since Kan extensions
along fully faithful functors are splittings of
restriction. Thus a map $F:\HO(Q)\to \HO(R)$ restricts
to a map $\rho(F):Q\to R$. In fact, we have a natural
isomorphism $\rho\circ\HO\iso\id_{\cat{QCAT}}$, so that
$\rho\circ\HO(f)$ is again $f$, up to this
isomorphism. Indeed, given $f:Q\to R$, we already know 
how to compute
$\HO(f)$ as $\Ho\circ \left(f^{N(-)}\right)$.
Then the restriction
$\rho(\HO(f)):Q\to R$, 
which we are to show coincides with $f$, is given
by $\rho(\HO(f))_n=\ob\circ \Ho\circ f^{\Delta^n}$.
That is, $\rho(\HO(f))$ acts by the action 
of $f$ on the objects of the homotopy categories 
of $Q^{\Delta^n}$ and $R^{\Delta^n}$. In other words,
it acts by the action of $f$ on the sets 
$\cat{SSET}(\Delta^n,Q)$ and $\cat{SSET}(\Delta^n,R)$;
via Yoneda, $\rho(\HO(f))$ acts by $f$ itself. 

It remains to show that $\HO(\rho(F))=F$
for any $F:\HO(Q)\to\HO(R)$.
By Lemma~\ref{catset lemma} it suffices to show that 
the restrictions of $\HO(\rho (F))$ and $F$ to the 
underlying $\cat{Dia}-$sets coincide. Using 
Lemma~\ref{kan lemma} and the adjunction characterizing the
Kan extension, we have
\[\cat{SET^{Dia^{\op}}}(\HO(Q)^{\ob}, \HO(R)^{\ob})
=\cat{SET^{Dia^{\op}}}(j_*Q,j_*R)\iso
\cat{SSET}(j^*j_*Q,R)\iso\cat{SSET}(Q,R)\]
In particular,
maps between $\HO(Q)^{\ob}$ and $\HO(R)^{\ob}$
agree when their restrictions to $Q$ and $R$ do.
Thus we are left to show that
$\rho(\HO(\rho(F)))=\rho(F).$
But as we showed above, $\rho\circ\HO$
is the identity map on $\cat{SSET}(Q,R)$,
so the proof is complete.
\end{proof}

\section{The embedding $\twocat{QCat}\to\twocat{PDer}$ of 2-categories}\label{section:2cat}

We shall now prove an analogous embedding theorem in the 2-categorical setting.
\begin{thm}\label{Mainthm2cat} 
Let $\twocat{QCat}$ denote the 2-category of small quasicategories.
Then 
the 2-functor $\HO:\twocat{QCat}\to\twocat{PDer}_{\twocat{Cat}}$ is 
bicategorically fully faithful; that is, it induces equivalences of hom-categories
$\twocat{QCat}(Q,R)\simeq \twocat{PDer}_{\twocat{Cat}}(\HO(Q),\HO(R))$ for any quasicategories
$Q$ and $R$. 	
\end{thm}

We get a Whitehead theorem for quasicategories as a corollary, following Remark 
\ref{ffconservative}. However, note that the following
is implied by Theorem \ref{whiteheadforquasicats}, whereas neither one of 
Theorem \ref{Mainthm2cat} and Theorem \ref{whiteheadforquasicats} implies the other.
\begin{cor}
The 2-functor $\HO:\twocat{QCat}\to \twocat{PDer}_{\twocat{Cat}}$ is bicategorically 
conservative; that is, if $f:Q\to R$ is a morphism of small quasicategories such that 
$\HO(f)$ is an equivalence in $\twocat{PDer}_{\twocat{Cat}}$, then $f$ is an equivalence
in $\twocat{QCat}$.
\end{cor}

The core tool for the proof is Theorem \ref{delocalization} below, 
which says that every quasicategory is a localization of a category. 
It is due to
Joyal but was first published by Stevenson in \cite{stevenson}.

First we recall the notion of $\infty$-localization, often just ``localization,"
for simplicial sets and quasicategories.
\begin{definition} Let $f:S\to T$ be a map of simplicial sets and 
$\mc W\subset S_1$ a set of edges. For any quasicategory $Q$, let
$Q^S_{\mc W}$ be the full sub-quasicategory of $Q^S$ on those maps
$g:S\to Q$ such that $g(w)$ is an equivalence in $Q$ for every edge
$w\in \mc W$. 

Then we say $f$ exhibits $T$ as an $\infty$-localization of $S$ at $\mc W$ if, for every 
quasicategory $Q$, the morphism $f^*:Q^T\to Q^S$ factors via an equivalence
$f^*:Q^T\to Q^S_{\mc W}$ of quasicategories.
\end{definition}

In particular, if $f:S\to T$ is a localization at $\mc W$ then for any quasicategory
$Q$, the pullback $f^*:\Ho(Q^T)\to\Ho(Q^S)$ is fully faithful, as we will use repeatedly
below. Specifically, $f^*$ is an equivalence onto the full subcategory 
$\Ho(Q^S_{\mc W})\subset \Ho(Q^S)$, since the 2-functor $\Ho$ preserves equivalences.

Let $\Delta\downarrow S$ be the category of simplices of a simplicial set $S$, and let
$p_S:N(\Delta~\downarrow~S)~\to~S$ be the natural extension of the projection 
$(f:\Delta^m\to S)\mapsto f(m)$. Finally, let $\mc L_S$ be 
the class of arrows $a:(f:\Delta^m\to S)\to (g:\Delta^n\to S)$ in $\ddown{S}$
such that $a(n)=m$, that is, the last-vertex maps. Then we have the following theorem:

\begin{thm}[\cite{stevenson}]\label{delocalization}
For any quasicategory $Q$, the last-vertex projection $p_Q$ exhibits $Q$ as an $\infty$-localization of the nerve
$N(\Delta\downarrow Q)$ at the class $\mc L_Q$.
\end{thm}

Thus every quasicategory $Q$ is canonically a localization of its category
$\Delta\downarrow Q$ of simplices. 

\begin{remark}
Observe that $N(\Delta\downarrow(-))$ constitutes an endofunctor of simplicial sets and that
$p:N(\Delta\downarrow(-))\to \mathrm{id}_{\cat{SSet}}$ is a natural transformation.
\end{remark}

We turn to the proof.

\begin{proof}[Proof of Theorem \ref{Mainthm2cat}]
First, we must show that if $F:\HO(Q)\to\HO(R)$ is a pseudonatural transformation,
then there exists $h:Q\to R$ and an isomorphism $\Lambda:\HO(h)\cong F$.
Observe that, since $Q$ is small, $\Delta\downarrow Q$ is in $\twocat{Cat}$.
Now we claim that $F_{\Delta\downarrow Q}(p_Q):$\mbox{$\ddown{Q}$}$\to R$ sends the
class $\mc L_Q$ of last-vertex maps into equivalences in $R$.
Indeed, if $\ell:\Delta^1\to \ddown{Q}$ is in $\mc L_Q$, then we have, using $F$'s respect for 2-morphisms and the
structure isomorphism $F_\ell$,

\[F_{[0]}(\dia(\ell^*p_Q))=\dia( F_{[1]}(\ell^*p_Q))\cong\dia(\ell^*F_{\Delta\downarrow J}(p_Q)).\]

Thus $\dia(\ell^*F_{\Delta\downarrow J}(p_Q))$ is an isomorphism in
$\Ho(R)$, since $\dia(\ell^*p_Q)$ is an isomorphism in $\Ho(Q)$.
Then using the delocalization theorem, we
can define $h:Q\to R$ as any map admitting an isomorphism 
$\s:h\circ p_Q\cong F_{\Delta\downarrow Q}(p_Q)$. 
From $\sigma$, we get an invertible
modification \[\HO(\sigma):\HO(h\circ p_Q)\Rightarrow \HO(F_{\ddown Q}(p_Q)):\HO(\ddown Q)\to \HO(R).\]

We now construct an invertible modification $\Lambda:\HO(h)\Rightarrow F:\HO(Q)\to\HO(R)$. Fixing 
$J\in \twocat{Cat}$ and $X:N(J)\to Q$, since $p_J^*:\HO(R)(J)\to \HO(R)(\ddown J)$ is fully faithful 
we can uniquely define $\Lambda_{X,J}:\HO(h)_J(X)\cong F_J(X)$ by giving
$p_J^*(\Lambda_{X,J})$. To wit, we require
$p_J^*(\Lambda_{X,J})$ to be the composition

\begin{align*}
p_J^*\HO(h)_J(X)&=h\circ X\circ p_J=h\circ p_Q\circ \ddown X \\&\iso F_{\ddown Q}(p_Q)\circ \ddown X
\iso F_{\Delta\downarrow J}(p_Q\circ \Delta\downarrow X)=F_{\Delta\downarrow J}(X\circ p_J)\\&\iso  F_J(X)\circ p_J.	
\end{align*}

The first isomorphism is a component of $\HO(\sigma)$, 
while the latter two are components of $F$. The naturality of $\Lambda_{J,X}$
in $X$ thus follows from the facts that $F$ is pseudonatural and that $\HO(\s)$ is a modification.
So we have natural isomorphisms $\Lambda_J:\HO(h)_J\Rightarrow F_J$ for each $J$. To verify that the
$\Lambda_J$ assemble into a modification, consider any 
$u:K\to J$. Then we must show that, for any $X:J\to Q$, the
diagram 
\[\begin{tikzcd} \HO(h)_J(X)\circ u\ar[d,equals]\ar[r,"\Lambda_{J,X}*u"]
&F_J(X)\circ u\ar[d,"F_u"]\\
\HO(h)_K(X\circ u)\ar[r,"\Lambda_{K,X\circ u}"]&F_K(X\circ u)\end{tikzcd}\]
commutes. Using, as always, full faithfulness of the pullback
along a localization, we may precompose with $p_K$. Then the
modification axiom is verified by the commutativity of the
following diagram: 

\[\begin{tikzcd}
     hXup_K\ar[d,equals]\ar[r,"\Lambda_{J,X}*up_K"]  &   F_J(X)up_K\ar[d,equals]\\
    hXp_J\ddown u\ar[d,equals]\ar[r,"\Lambda_{J,X}*p_J\ddown u"]&   F_J(X) p_J\ddown u  \ar[r,equals] &   F_J(X)up_K\ar[from=d,"F_u*p_K"]\\
    hp_Q\ddown Xu\ar[d,equals]&     F_{\ddown J}(Xp_J)\ddown u\ar[u,"F_{p_J}*\ddown u"]\ar[d,equals]&    F_K(Xu)p_K\ar[from=d,"F_{p_K}"]\\
    F_{\ddown Q}(p_Q)\ddown{Xu}\ar[r,"F^{-1}_{\ddown X}*\ddown u"]& F_{\ddown J}(p_Q\ddown X)\ddown u&F_{\ddown K}(Xup_K)\\
    &F_{\ddown K}(p_Q\ddown{Xu})\ar[from=ul,"F^{-1}_{\ddown{Xu}}"]\ar[u,"F_{\ddown u}"]\ar[ur,equals]
\end{tikzcd}\]
The upper left square commutes since $up_K=p_J\ddown{u}$. The left central hexagon commutes by definition of 
$\Lambda_{J,X}$, and the lower left triangle and right-hand heptagon commute by functoriality of the 
pseudonaturality isomorphisms of $F$. Meanwhile, the outer route aruond the diagram from $hXup_K$ to 
$F_J(X)up_K$ is $F_u\Lambda_{K,Xu}$, while the inner route is $\Lambda_{J,X}*up_K$. So $\Lambda$ is an
invertible modification $\HO(h)\iso F$, as desired.

We have shown that $\HO$ induces an essentially surjective functor 
$\twocat{QCat}(Q,R)\to \twocat{PDer}(\HO(Q),\HO(R))$. We next consider full faithfulness.
 So, assume given a modification 
 \[\Xi:\HO(f)\Rightarrow\HO(g):\HO(Q)\to\HO(R).\]
 We must show there exists a unique $\xi:f\Rightarrow g$ with $\HO(\xi)=\Xi$.  
First, we consider
\[\Xi_{p_Q}:f\circ p_Q\to g\circ p_Q,\] 
which is a morphism in $\HO(R)(\ddown Q).$ 
According to \ref{Der5'},
we can lift this to a map 
$\widehat{\Xi}_{p_Q}:\Delta\downarrow Q\to R^{\Delta^1}$ with 
$\dia(\widehat{\Xi}_{p_Q})=\Xi_{p_Q}.$ 

Since the domain and codomain $f\circ p_Q$ and $g\circ p_Q$ of 
$\widehat{\Xi}_{p_Q}$ invert the last-vertex maps $\mc L_Q$, by (Der\ref{Der2}) so does 
$\widehat{\Xi}_{p_Q}$ itself. Thus by the delocalization theorem
we get 
$\widehat{\Xi}':Q\to R^{\Delta^1}$ with an isomorphism 
\[a:\widehat\Xi'\circ p_Q\iso \widehat\Xi_{p_Q}.\]
The domain and codomain 
\[0^*a:0^*(\widehat\Xi'\circ p_Q)\iso fp_Q \text{ and } 1^*a:1^*(\widehat\Xi'\circ p_Q)\iso gp_Q\] 
give rise to 
unique isomorphisms 
\[i:0^*\widehat\Xi'\iso f \text{ and }j:1^*\widehat \Xi'\iso g.\] 
Now we can construct $\widehat\Xi:Q\to R^{\Delta^1}$ 
as a lift of the composite 
\[\begin{tikzcd}f\ar[r,"i^{-1}"]&0^*\widehat \Xi'\ar[r,"\dia(\widehat \Xi')"]&1^*\widehat\Xi'\ar[r,"j"]&g\end{tikzcd}\] in 
$\Ho(R^Q)$. Using the fullness clause of \ref{Der5'}, we can choose an isomorphism 
$b:\widehat\Xi\iso\widehat\Xi'$ in $\Ho((R^{\Delta^1})^Q)$
lifting $(i^{-1},j^{-1}):\dia(\widehat\Xi)\to\dia(\widehat\Xi')$. 

Then 
$a\circ(b*p_Q):\widehat\Xi\circ p_Q\to \widehat{\Xi}_{p_Q}$ is an isomorphism with endpoints fixed, insofar as
$0^*(b*p_Q)=i^{-1}*p_Q=0^*a^{-1}$ and similarly $1^*(b*p_Q)=1^*a^{-1}$. Thus $\dia(\widehat\Xi\circ p_Q)=\dia(\widehat \Xi_{p_Q})=\Xi_{p_Q}$ in $\Ho(R^{\ddown Q})$. 

Notice that if $\widehat\Xi_2:Q\to R^{\Delta^1}$ is any other morphism satisfying 
$\dia(\widehat\Xi_2\circ p_Q)=\Xi_{p_Q}$, then $\dia(\widehat\Xi_2)=\dia(\widehat\Xi)$, since 
pullback along $p_Q$ is faithful. So we have a unique candidate $\xi:=\dia(\widehat\Xi):f\Rightarrow g$;
it remains to show that $\HO(\xi)=\Xi$. 

To that end, we claim that for every $X:J\to Q$, we have $\HO(\xi)_X=\xi * X=\Xi_X$. 
As above, it suffices to precompose $X$ with $p_J$, and then we have
\[\xi*X*p_J=\dia(\widehat\Xi)*p_Q*\ddown X=
\dia(\widehat\Xi\circ p_Q)\circ \Delta\downarrow X\]
\[=\Xi_{p_Q}*\Delta\downarrow X=\Xi_{p_Q\circ \ddown X}=\Xi_{X\circ p_J}=\Xi_X*p_J\]
as desired. In the equations above we have used the 2-functoriality of $\HO(R)$, naturality of $p$,
and the modification property of $\Xi$. So $\HO(\xi)=\Xi$, as was to be shown.
\end{proof}


\section{Whitehead's theorem for quasicategories}\label{section:whitehead}

As discussed above, Theorem \ref{Mainthm2cat} says that $\HO$ is bicategorically fully faithful,
and in particular bicategorically conservative, when the domain is small quasicategories and
$\twocat{Dia}=\twocat{Cat}$. In this section, we show that $\HO$ is 
at least bicategorically conservative no matter what assumptions are placed on the domain and
codomain. We first give informal arguments that this is the best one should hope for.

\subsection*{Evidence against a stronger result}

We do not hope to prove all of Theorem \ref{Mainthm2cat} for arbitrary 
choices of domain and codomain, as is most intuitive to see in the case of 
$\HO:\twocat{QCat}\to\twocat{PDer}_{\twocat{HFin}}$. 
Since the 2-category $\twocat{HFin}$ of homotopy finite categories is small,
prederivators of that domain form a ``concrete 2-category"
in the strongest possible sense. That is, we have a 2-functor $U:\twocat{PDer}_{\twocat{HFin}}\to \twocat{Cat}$, faithful on 1- and 2-morphisms, given by 
\[U(\D)= \prod_{J\in\twocat{HFin}}\D(J)\times\prod_{u:K\to J}\D(J)^{[1]}\]
For $F:\D_1\to\D_2$, we have $U(F)=((F_J),(F_u:\D(J)^{[1]}\to\D(K)^{[1]}))$, 
where the functor $F_u$ sends $f:X\to Y$ to the arrow $u^*F(X)\to F(u^*Y)$ 
which can be defined in two equivalent ways using the pseudonaturality
isomorphisms of $F$.
Similarly, for $\Xi:F_1\Rightarrow F_2$, we have 
$U(\Xi)=((\Xi_J),(\Xi_u:F_u\Rightarrow G_u))$, where the components of $\Xi_u$ 
are $u*\Xi_X$ and $\Xi_{u^*Y}$. It is straightforward to check
that these objects are, respectively, a functor and a natural transformation. 
\footnote{The reason for the unfamiliar $u$ terms in the definition of $C$ is
that a pseudonatural transformation is not determined by its action
on objects.}

Since $\HO:\twocat{QCat}\to\twocat{PDer}_{\twocat{HFin}}$ is faithful
on 1-morphisms, if it were also faithful on 2-morphisms then $\twocat{QCat}$ 
would be a strictly concrete 2-category in the same sense as $\twocat{PDer}_{\twocat{HFin}}$. In perhaps more familiar terms, there would be 
no ``phantom homotopies" between maps of quasicategories. 
That this should be the case strains credulity, given the famous
theorem of Freyd \cite{freyd} that the category of spaces $\cat{Hot}$ is not
concrete (i.e. admits no faithful functor to $\cat{Set}$,) 
though we propose no specific counterexample.

We do conjecture a specific counterexample to the claim that 
$\HO:\twocat{QCat}\to \twocat{PDer}_{\twocat{HFin}}$ is locally 
essentially surjective. That is, we suggest a pseudonatural transformation 
$F:\HO(Q)\to \HO(R)$ not isomorphic to the image of any quasicategory morphism
$f:Q\to R$. A mapping telescope argument rules out any countable $Q$, so we turn to the
simplest possible uncountable example.

Denote then, as usual, the least uncountable
ordinal by $\omega_1$. We construct a ``partially coherent" form of $\omega_1$ 
as follows. We have a functor $D:\omega_1\to \cat{SSet}$ 
that sends each countable ordinal $\a\in\omega_1$ to its nerve and each
map $\a<\b$ in $\omega_1$ to the standard inclusion $N\a\mono N\b$. Since
$D$ is projectively cofibrant, the homotopy
colimit of $D$ in $\cat{SSet}_{\mathrm{Joyal}}$ is equivalent both to 
$N\omega_1$ and to the homotopy colimit of the simplicial set 
$n\mapsto \coprod_{\a_0<...<\a_n} N(\a_n)$. 
We define $\omega_1^{A_2}$ as the homotopy colimit in 
$\cat{SSet}_{\mathrm{Joyal}}$ of the restriction of this simplicial object to 
$\Delta^\op_{\leq 2}$, the full subcategory of $\Delta^\op$ on the objects $[0],[1],[2]$.

Thus, concretely, $\omega_1^{A_2}$ is given
by first taking the disjoint union of the nerves of the countable ordinals 
$\a$, then gluing in an isomorphism from $\a$ to its image in $\b$ whenever 
$\a<\b$, and finally gluing in a commutative triangle between these isomorphisms for each triple 
$\a<\b<\g$. A map from $\omega_1^{A_2}$ into a quasicategory thus corresponds to 
maps out of each countable ordinal which are homotopic via homotopies 
which commute, but which need not satisfy any higher coherences. This explains the
notation $\omega_1^{A_2}$ as an analogy between a (fully coherent) homotopy 
colimit and an $A_\infty$-algebra. 

There is a canonical pseudonatural transformation 
$F:\HO(\omega_1)\to \HO(\omega_1^{A_2})$ defined as follows. Given any 
homotopically finite category $J$ and a map $X:N(J)\to N(\omega_1)$, 
let $\a\in\omega_1$ be the least ordinal through which $X$ factors. Then
we define $F(X):N(J)\to \omega_1^{A_2}$ as the composition of $X$ with the
canonical inclusion of $\a$ in $\omega_1^{A_2}$. The pseudonaturality 
constraints of $F$ are constructed from the canonical isomorphisms in
$\omega_1^{A_2}$. 

\begin{conj}\label{conj:incoherent}
The pseudonatural transformation $F$ defined above is not isomorphic to a
strict transformation, and is thus not in the essential image of $\HO$. 
\end{conj}

To prove the conjecture, it would suffice to show that every
morphism $N(\omega_1)\to \omega_1^{A_2}$ is bounded in terms of the
countable ordinals $\a$ whose image in $\omega_1^{A_2}$ it intersects. 
This boundedness should follow from analysis of the Kan complex under 
$\omega_1^{A_2}$ produced by collapsing the image of each $\a$ to a point. 

\subsection*{Proof of Whitehead's theorem for quasicategories}
We will use the main 
theorem of 
\cite{whiteheadforspaces}, which says that the 2-category $\twocat{KAN}\subset\twocat{QCAT}$ of Kan 
complexes is strongly generated by the tori $(S^1)^n$, in the sense that a morphism $f:X\to Y$ of 
Kan complexes is a homotopy equivalence if and only if, for each $n$, the functor 
$\twocat{KAN}((S^1)^n,f)$ is an equivalence of groupoids.\footnote{We write $S^1$ for any Kan complex of the homotopy type of the circle.} We rephrase this in a form more convenient
for our purposes:

\begin{thm}\label{whiteheadforspaces}
The restriction of $\HO:\twocat{QCAT}\to\twocat{PDer}_{\twocat{HFin}}$ to the 2-category 
$\twocat{KAN}$ reflects equivalences.
\end{thm}

Recall that equivalences in $\twocat{PDer}_{\twocat{HFin}}$ in the abstract 2-categorical sense
coincide with pseudonatural transformations which induce equivalences of categories levelwise.

\begin{proof}
Given $f:X\to Y$ in $\twocat{KAN}$, the image $\HO(f)$ is an equivalence in $\twocat{PDer}_{\twocat{HFin}}$
if and only if, for every homotopically finite category $J$, the induced functor $\Ho(f^{N(J)}):\Ho(X^{N(J)})\to \Ho(Y^{N(J)})$ is an equivalence. Since the classical model structure on simplicial sets is also
Cartesian, we have equivalences $\Ho(X^{N(J)})\equiv \Ho(X^{\mathrm{Ex}^\infty(N(J))})$, 
and similarly for $Y$,
where $\mathrm{Ex}^\infty$ is Kan's fibrant replacement functor. Now, by Thomason's theorem \cite{thomason}, 
as $J$ varies, $\mathrm{Ex}^\infty(N(J))$ runs through all finite homotopy types. In particular, if $\HO(f)$ is an
equivalence in $\twocat{PDer}$, then $f$ induces equivalences 
$\Ho(X^{(S^1)^n})\to\Ho(Y{(S^1)^n})$ 
for every $n$, which is to say, $\twocat{KAN}((S^1)^n,f)$ is an equivalence. Thus $f$ must be an 
equivalence, by \cite{whiteheadforspaces}.
\end{proof}

To make use of the above result to prove results on the relationship between 
quasicategories and their prederivators, we first recall what Rezk has 
described as the fundamental theorem of quasicategory theory. 
First, a quasicategory $Q$ has mapping spaces $Q(x,y)$ for each $x,y\in Q$,
which can be given various models. We shall use the balanced model in which
we have $Q(x,y)=\{(x,y)\}\times_{Q\times Q} Q^{\Delta^1}$, so that an 
$n$-simplex of $Q(x,y)$ is a prism $\Delta^n\times\Delta^1$ in $Q$ which
is degenerate on $x$ and $y$ at its respective endpoints.

We say that a map $f:Q\to R$ of quasicategories is \emph{fully faithful} if 
it induces an equivalence of Kan complexes $Q(x,y)\to R(f(x),f(y))$ for 
every $x,y\in Q$. It is essentially surjective if, for every $z\in R$,
there exists $x\in Q$ and an edge $a:f(x)\to z$ which becomes an isomorphism 
in $\Ho(R)$. Then we have
\begin{thm}[Joyal]
A map $f:Q\to R$ of quasicategories is an equivalence in the sense of
Definition \ref{equivalenceofquasicategories} if and only if it is
fully faithful and essentially surjective.	
\end{thm}

Now we can prove our Whitehead theorem for quasicategories. 

\begin{thm}\label{whiteheadforquasicats}	
 Let $f:Q\to R$ be a map of quasicategories, and suppose that $\HO(f)$ is
an equivalence in $\twocat{PDer}$. Then $f$ is an equivalence of quasicategories.

\end{thm}

\begin{proof}
Since $\Ho(f)$ is an equivalence by assumption, $f$ is essentially surjective.
Thus we have only to show $f$ is fully faithful. By Theorem 
\ref{whiteheadforspaces},
it suffices to show that $\HO(f)$ induces an equivalence
$\HO(Q(x,y))\cong \HO(R(f(x),f(y)))$ in $\twocat{PDer}$ for every $x$ and $y$ in $Q$. 
What is more, since for any $J$ we have $Q(x,y)^{NJ}\iso Q^{NJ}(p_J^*x,p_J^*y)$, 
it suffices at last to show that $f$ induces equivalences 
$f_{x,y}:\Ho(Q(x,y))\to \Ho(R(f(x),f(y)))$ on the homotopy categories of 
mapping spaces.

Essential surjectivity is proved via an argument that also 
appeared in the construction of $\widehat\Xi$ in the proof of Theorem
\ref{Mainthm2cat}. Namely, from essential surjectivity of $\HO(f)$, 
given any $X\in \Ho(R(f(x),f(y)))$ and any $Y\in \HO(Q)([1])$ with an 
isomorphism $s:\HO(f)(Y)\cong X$ in $\HO(R)([1])$, we see by 
conservativity and fullness of $\HO(f)$ that we have isomorphisms
$0^*Y\iso x\in\HO(Q)(J)$ and, similarly, $1^*Y\iso y$. Composing these 
isomorphisms and $\dia Y$ in $\Ho(Q)$ gives a morphism $x\to y$ in $\Ho(Q)$ 
isomorphic to $\dia Y$ in $\Ho(Q)^{[1]}$.
By (Der5') and (Der\ref{Der2}) we can lift this to an isomorphism 
$r:Y'\cong Y$ in $\HO(Q)([1])$ such that $0^*(s\circ\HO(f)(r))=\id_{x}$ and
$1^*(s\circ\HO(f)(r))=\id_{y}$. This implies that $s\circ \HO(f)(r)$
may be lifted to an isomorphism $\HO(f)(Y')\iso X$ in $\Ho(R(f(x),f(y))$.
Thus $f_{x,y}$ is essentially surjective. 

For fullness, we observe that if 
$a:Y_1\to Y_2\in \HO(Q)([1])$ verifies $Y_1,Y_2:x\to y$, 
$0^*\HO(f)(a)=\id_{f(x)}$, and $1^*\HO(f)(a)=\id_{f(y)}$, then 
we have also $0^*(a)=\id_x$ and $1^*a=\id_y$, since $\HO(f)$ is faithful. 
This implies that $a$ can be lifted to a morphism $a':Y_1\to Y_2$ in 
$\Ho(Q(x,y))$ with $f_{x,y}(a)=\HO(f)(a)$. And since $\HO(f)$ is full, every
morphism $\HO(f)(Y_1)\to \HO(f)(Y_2)$ in $\Ho(R(f(x),f(y)))$ is equal to
$\HO(f)(a)$ in $\HO(R)([1])$, for some $a$. 

Finally, we turn to faithfulness. 
Suppose we have morphisms $a,b:Y_1\to Y_2$ in $\Ho(Q(x,y))$
with $f_{x,y}(a)=f_{x,y}(b)$ in $\Ho(R(f(x),f(y)))$. 
We wish to show $a=b$. First, we may represent $a$ and $b$ by 
$\hat a,\hat b\in \HO(Q)( [1]\times [1])$, each with boundary 
\[
\begin{tikzcd}
x\ar[r,"Y_1"]\ar[d,equals]&y\ar[d,equals]\\
x\ar[r,"Y_2"]&y
\end{tikzcd}
\]

Let $\del[2]$ denote  the category on objects $0,1,2$ freely generated by 
\emph{three} arrows $0\to 1,1\to 2,0\to 2$, so that $N\del[2]$ is Joyal 
equivalent to $\del\Delta^2$. 
The lifts $\hat a$ and $\hat b$ fit together in a 
diagram $W\in\HO(Q)([1]\times \del [2])$ with 
$(01)^* W=q^*Y_1$, $(02)^*W=\hat a$, and 
$(12)^*W=\hat b$, where $q:[1]\times[1]\to [1]$ projects out
the last coordinate.  The significance of $W$ is that we have $a=b$ if and only
if $W$ admits an extension $Z$ to $\HO(Q)([1]\times [2])$ such that
$Z|_{\{0\}\times [2]}=p_{[2]}^* x$ and $Z|_{\{1\}\times [2]}=p_{[2]}^* y$. 
It suffices to exhibit $W'\in \HO(Q)([1]\times \del[2])$ with 
$W'|_{0\times \del[2]}=p_{\del [2]}^*x$ and 
$W'|_{1\times \del[2]}=p_{\del [2]}^*y$ admitting such
an extension $Z'$, together with an isomorphism 
$t:W\to W'$ in $\HO(Q)([1]\times\del [2])$ 
such that $t|_{{0}\times \del[2]}=\id_{p_{\del [2]}^* x}$ and
$t|_{{1}\times \del[2]}=\id_{p_{\del [2]}^* y}$. Indeed, in this situation 
$W$ and $W'$ both represent maps from $S^1$ to the Kan complex 
$Q(x,y)$, $Z$ and $Z'$ represent putative extensions to $\Delta^2$, 
and $t$ represents a homotopy between them. 

In particular, since by assumption $\HO(f)(a)=\HO(f)(b)$ in 
$\Ho(R(f(x),f(y)))$, there exists an extension
$T$ of $\HO(f)(W)$ to $\HO(R)([1]\times [2])$ with trivial endpoints,
as above. Now take $\hat T\in \HO(Q)( [1]\times[2])$ with 
an isomorphism $s:\HO(f)(\hat T)\iso T$. In particular, this gives isomorphisms 
$\HO(f)(\hat T)|_{\{0\}\times [2]}\iso p_{[2]}^*f(x)$ 
and $\HO(f)(\hat T)|_{\{1\}\times [2]}\iso p_{[2]}^*f(y)$ in $\HO(R)([2])$, 
which lift uniquely to isomorphisms $\hat T|_{\{0\}\times [2]}\iso p_{[2]}^*x$ 
and $\hat T|_{\{1\}\times [2]}\iso p_{[2]}^*y$ in $\HO(Q)([2])$. 
Composing these isomorphisms
with $\dia\hat T$ and lifting into $\HO(Q)([1]\times [2])$ 
gives $Z'\in \HO(Q)([1]\times [2])$ with $Z'|_{\{0\}\times [2]}=p_{[2]}^*x$ and
$Z'|_{\{1\}\times [2]}=p_{[2]}^*y$, together with an isomorphism 
$t':\HO(f)( Z')\iso T$ in $\HO(R)([1]\times[2])$ inducing the identity on 
$p_{[2]}^*f(x)$ and $p_{[2]}^*f(y)$, respectively. Restricting $t'$ to 
$[1]\times \del[2]$ and lifting to $\HO(Q)([1]\times\del[2])$ 
specifies an isomorphism $t:Z'|_{[1]\times\del[2]}\iso W$ such that 
$t|_{{0}\times \del[2]}=\id_{p_{\del [2]}^* x}$ and
$t|_{{1}\times \del[2]}=\id_{p_{\del [2]}^* y}$. As we saw above, this suffices
to guarantee that $W$ admits an extension $Z$ as desired. 
\end{proof}

\bibliography{bib}{}
\bibliographystyle{plain}

\end{document}